\documentclass[11pt]{amsart}

\usepackage{amsmath,amsfonts,amssymb,amscd}
\begin{document}
\renewcommand{\thefootnote}{\fnsymbol{footnote}}
\pagestyle{plain}

\title{A remark on the Donaldson-Futaki invariant \\
for sequences of test configurations}
\author{Toshiki Mabuchi${}^*$}
\maketitle
\footnotetext{ ${}^{*}$Supported 
by JSPS Grant-in-Aid for Scientific Research (A) No. 20244005.}
\abstract
In this note, we consider a sequence of test configurations $\mu_j = (\mathcal{X}_j,\mathcal{L}_j,\psi_j )$, 
$j=1,2,\dots$,   
compatible with a K\"ahler metric $\omega = c_1(L;h)$ on a polarized algebraic manifold $(X,L)$.
Then an explicit formula for the Donaldson-Futaki invariant $F_1(\{\mu_j\})$ in \cite{M} 
will be given.
 
\endabstract
\section{Introduction}

By a {\it polarized algebraic manifold} $(X,L)$, we mean
a pair of a nonsingular irreducible projective variety $X$, 
defined over $\Bbb C$, and a very ample line bundle $L$
over $X$. 
In this note, we fix once for all such a pair $(X,L)$ with $n:= \dim X$. For the affine line 
 $ \Bbb A^1 := \{z\in \Bbb C\}$, we consider the one-dimensional algebraic torus $\Bbb G_m :=\{t\in \Bbb C^*\}$ acting on $\Bbb A^1$ by multiplication of complex numbers
$$
\Bbb G_m\times \Bbb A^1 \to \Bbb A^1,
\qquad (t, z) \mapsto
t z.
$$
Fix a Hermitian metric $h$ for $L$ such that $\omega := c_1(L;h)$ is K\"ahler.
Endow $V_{\ell}:= H^0(X,L^{\otimes \ell})$, $\ell =1,2,\dots$, 
with a Hermitian metric $\rho_{\ell}$ defined by
$$
\langle \sigma', \sigma''\rangle_{\rho_{\ell}}\; :=\; \int_X \,(\sigma', \sigma'')_h \,\omega^n,
\qquad \sigma', \sigma'' \in V_{\ell},
\leqno{(1.1)}
$$
where $(\sigma', \sigma'')_h$ denotes the pointwise Hermitian inner product of $\sigma'$ and $\sigma''$ 
by the $\ell$-multiple of $h$.
For the Kodaira embedding 
$$
 \Phi_{\ell} \,:\, X \, \hookrightarrow \, \Bbb P^*(V_{\ell})
 $$
associated to the complete linear system $|L^{\otimes \ell}|$ on $X$,
we consider its image $X_{\ell}:= \Phi_{\ell}(X)$, and choose an algebraic group homomorphism
$$
\psi \, : \,\Bbb G_m \,\to \,\operatorname{GL}(V_{\ell})
$$
such that the maximal compact subgroup $S^1 \subset \Bbb G_m$
acts isometrically on $(V_{\ell}, \rho_{\ell})$.
Let $\mathcal{X}^{\psi}$ be the irreducible algebraic subvariety of $\Bbb A^1 \times \Bbb P^* (V_{\ell})$ obtained as 
the closure of $\cup_{z\in \Bbb G_m} \mathcal{X}^{\psi}_z$ in $\Bbb A^1 \times \Bbb P^* (V_{\ell})$ 
by setting
$$
\mathcal{X}^{\psi}_z := \{z\}\times\psi (z)   X_{\ell},
\qquad z \in \Bbb G_m,
$$
where the element $\psi (z)  $ in $\operatorname{GL}(V_{\ell})$ 
acts naturally on  the 
set $\Bbb P^* (V_{\ell})$ of all hyperplanes in $V_{\ell}$ passing through the origin.
We then consider the map
$$
\pi : \mathcal{X}^{\psi} \to \Bbb A^1
$$ 
induced by the projection  $\operatorname{pr}_1 : \Bbb A^1 \times \Bbb P^* (V_{\ell}) \to \Bbb A^1$ 
to the first factor.
For the hyperplane bundle $\mathcal{O}_{\Bbb P^*(V_{\ell})}(1)$ on $\Bbb P^*(V_{\ell})$, we consider the restriction 
$$
\mathcal{L}^{\psi}\, :=\,\operatorname{pr}_2^*\mathcal{O}_{\Bbb P^*(V_{\ell})}(1)_{|\mathcal{X}^{\psi}},
$$
where $\operatorname{pr}_2 : \Bbb A^1 \times \Bbb P^* (V_{\ell}) \to \Bbb P^* (V_{\ell})$
denotes the projection to the second factor.
For the dual vector space $V_{\ell}^*$ of $V_{\ell}$,
the $\Bbb G_m$-action on $\Bbb A^1 \times V_{\ell}^*$ defined by
$$
\Bbb G_m \times (\Bbb A^1 \times V_{\ell}^*)\to \Bbb A^1 \times V_{\ell}^*,
\quad (t, (z, p))\mapsto  (tz, \psi (t) p),
$$
naturally induces $\Bbb G_m$-actions on $\Bbb A^1 \times \Bbb P^*(V_{\ell})$ and $\mathcal{O}_{\Bbb P^*(V_{\ell})}(-1)$, where $\operatorname{GL}(V_{\ell})$ 
acts  on $V_{\ell}^*$ 
by contragradient representation.
 This then induces $\Bbb G_m$-actions on $\mathcal{X}^{\psi}$ and $\mathcal{L}^{\psi}$, 
while $\pi : \mathcal{X}^{\psi} \to \Bbb A^1$ above is a projective morphism with relative very ample line bundle 
$\mathcal{L}^{\psi}$ such that
$$
(\mathcal{X}^{\psi}_z, \mathcal{L}_z^{\psi})\; \cong \; (X,L^{\otimes \ell}),
\qquad z \neq 0,
$$
where for each $z \in \Bbb A^1$, we denote by $\mathcal{L}_z^{\psi}$ the restriction of $\mathcal{L}^{\psi}$ to 
the scheme-theoretic fiber $\mathcal{X}^{\psi}_z := \pi^{-1}(z)$.
Then a  triple $\mu = ({\mathcal{X}}, {\mathcal{L}}, \psi )$ is called 
a {\it
test configuration for $(X,L)$}, if we have both 
$$
\mathcal{X} = \mathcal{X}^{\psi}\quad\text{ and  }\quad\mathcal{L}= \mathcal{L}^{\psi}
$$
for some $\psi$ as above,
where $\ell$  is called 
the {\it exponent} of the test configuration $\mu = ({\mathcal{X}}, {\mathcal{L}}, \psi )$.
For positive integers $\gamma$, let $k = \gamma \ell$.
Then for the homogeneous ideal  $I = \oplus_{\gamma} \, I_{\gamma}$ for $\mathcal{X}_0$ in $\Bbb P^*(V_{\ell})$, 
we put
$$
R_k \;:=\; S^{\gamma}(V_{\ell})/I_{\gamma}, \;\;\qquad \gamma =1,2,\dots,
$$
where $S^{\gamma}(V_{\ell})$ denote the $\gamma$-th symmetric tensor product of $V_{\ell}$.
Since the $\Bbb G_m$-action on $V_{\ell}$ preserves $I_{\gamma}$, we have a natural $\Bbb G_m$-action 
on $R_k$. Let $w_{k}(\mathcal{X},\mathcal{L})$ be the weight of the $\Bbb G_m$-action on $R_k$. Put $N_{k} := \dim R_k$. 
Then
$$
F_0 (\mathcal{X},\mathcal{L} )\; := \; \lim_{k\to \infty} \,\frac{w_k(\mathcal{X},\mathcal{L})}{k N_k}\; =\; \lim_{k \to \infty}
\frac{n!\, w_k (\mathcal{X},\mathcal{L})}{k^n\, c_1(L)^n[X]}
$$
is the $F_0$-term in the Donaldson-Futaki asymptotic expansion for the {\it test configuration
$(\mathcal{X},\mathcal{L})$ in Donaldson's sense}.
For $\gamma =1$, 
$N_{\ell}$ coincides with $\dim V_{\ell}$, while
$w_{\ell}(\mathcal{X},\mathcal{L})$ is the weight of the $\Bbb G_m$-action on $V_{\ell}$.
For the real Lie subgroup $\Bbb R_+ := \{  t \in \Bbb R_+\}$ of $\Bbb G_m$, 
we define an associated Lie group homomorphism $\psi^{\operatorname{SL}}: \Bbb R_+ \to 
\operatorname{SL}(V_{\ell})$ by
$$
\psi^{\operatorname{SL}}(t) \;:=\; \frac{\psi (t)}{\det ( \psi (t) )^{1/N_{\ell}}}, 
\qquad t \in \Bbb R_+.
$$
As in \cite{M}, for a test configuration $\mu = (\mathcal{X},\mathcal{L},\psi )$ for $(X,L)$ of exponent $\ell$, 
we define  $\|\mu \|_{\infty}$ and $\|\mu \|_1$ by
$$
\begin{cases}
&\|\mu \|_1 \, := \;(|b_{1}| + |b_2| + \dots + |b_{N_{\ell}}|)/\ell^{n+1},\\
&\|\mu \|_{\infty}\, :=\;\max \{|b_{1}|, |b_2|, \dots , |b_{N_{\ell}}|\}/\ell,
\end{cases}
$$
where $b_{\alpha}$, $\alpha = 1,2,\dots, N_{\ell}$, are the weights of the $\Bbb R_+$-action 
via $\psi^{\operatorname{SL}}$ on the dual vector space $V^*_{\ell}$.
Put $\lambda := c_1(L)^n[X]/\pi$.
Consider 
the set $\mathcal{M}$ of all sequences $\{\mu_j \}$ of 
test configurations 
$$
\mu_j  \, =\, (\mathcal{X}_j, \mathcal{L}_j, \psi_j ), \qquad j =1,2,\dots,
\leqno{(1.2)}
$$
for $(X,L)$ such that the exponent $\ell_j$ of the 
test configuration $\mu_j$ satisfies the following 
growth condition:
$$
\text{$\ell_j \to +\infty$, \; as $j \to \infty$.}
$$
Here for simplicity, we assume that all $\mu_j$ are nontrivial. 
In \cite{M}, we defined the Donaldson-Futaki invariant $F_1 (\{\mu_j\})$ for sequences 
$\{\mu_j\}$ of test configurations. 
The purpose of this note is to show the following:
$$
F_1 (\{\mu_j \} ) \; =\; \varliminf_{j\to \infty}\, 
\frac{\lambda}{\|\mu_j\|_1}\left \{ \frac{w_{\ell_j}(\mathcal{X}_j,\mathcal{L}_j)}{N_{\ell_j}} - \ell_j F_0(\mathcal{X}_j,\mathcal{L}_j )\right \}.\quad
\leqno{\bf Theorem.}
$$
In the situation  as in Fact of \cite{M}, 
we have $\mathcal{X}_j = \mathcal{X}$ and $\mathcal{L}_j = \mathcal{L}^{\otimes \ell_j}$.
Hence $F_0(\mathcal{X}_j,\mathcal{L}_j ) = F_0 (\mathcal{X}, \mathcal{L})$  
and $w^{}_{\ell_j} (\mathcal{X}_j,\mathcal{L}_j ) = w^{}_{\ell_j} (\mathcal{X}, \mathcal{L})$ for all $j$.
Then, since
\begin{align*}
&\lim_{j\to \infty} \left\{ \frac{w_{\ell_j}(\mathcal{X}_j,\mathcal{L}_j)}{N_{\ell_j}} -\ell_j F_0 (\mathcal{X}_j,\mathcal{L}_j)\right \}\\
&= \lim_{j\to \infty} \ell_j\left \{  \frac{w_{\ell_j}(\mathcal{X},\mathcal{L})}{\ell_j N_{\ell_j}}
  - F_0(\mathcal{X},\mathcal{L})\right \} =
 F_1 (\mathcal{X},\mathcal{L}),
\end{align*}
and since $\beta = \lim_{j\to \infty} \|\mu_j\|_1$, we obtain (1.5) in \cite{M} from Theorem above. 
Later in (4.15) and Proposition of Section 4, we shall obtain a reformulation of Theorem E in 
Appendix of \cite{M4}.

\section{The Donaldson-Futaki invariant for sequences}

For a sequence $\mu_j = (\mathcal{X}_j, \mathcal{L}_j, \psi_j)$, $ j = 1,2,\dots$, of test configurations
of exponent $\ell_j$ as in (1.1), we
put $d_j := \ell_j^n\,c_1(L)^n[X]$ and
$W_j := \{S^{d_j}(V_{\ell_j})\}^{\otimes n+1}$. Then the dual space $W^*_j$ of $W_j$ 
admits the Chow norm (see \cite{Zh})
$$
W_j^* \owns w \,\mapsto \, \|w\|_{\operatorname{CH}(\rho^{}_{\ell_j})} \in \Bbb R_{\geq 0},
$$
associated to the Hermitian metric $\rho^{}_{\ell_j}$ for $V_{\ell_j}$.
For the Kodaira embedding $\Phi_{\ell_j}: X \hookrightarrow \Bbb P^*(V_{\ell_j})$ as in the introduction, let 
$$
0 \neq \hat{X}_{\ell_j}\in W_j^*
$$
be the associated Chow form for $X_{\ell_j} := \Phi_{\ell_j}(X)$ viewed as an irreducible reduced 
algebraic cycle on the projective space $\Bbb P^*(V_{\ell_j})$.
Let $\delta_j$ be such that
$$
\delta_j \; :=\; 
\begin{cases}
&\|\mu_j\|_{\infty}/\|\mu_j\|_1, \quad \text{ if $\|\mu_j\|_{\infty} \neq 0$;} \\
&  \;\; 1 \qquad \qquad\qquad  \,\text{ if $\|\mu_j\|_{\infty} = 0$;}
\end{cases}
$$
If $\|\mu_j\|_{\infty} \neq 0$, each $t$ in the real Lie subgroup $\Bbb R_+$ of $\Bbb G_m $ 
will be reparametrized by $t = \exp (s/\| \mu_j \|_{\infty})$ 
for some $s\in \Bbb R$, 
while if $\|\mu_j\|_{\infty} = 0$, no relations between $s\in \Bbb R$ and $t\in \Bbb G_m$ are required. Note that
$$
\psi_{j}^{\operatorname{SL}}: \Bbb R^+ \to \operatorname{SL}(V_{\ell_j})
$$
is the restriction of $\psi_j : \Bbb G_m \to \operatorname{SL}(V_{\ell_j})$ to $\Bbb R_+$.
Since the group $\operatorname{SL}(V_{\ell_j})$ 
acts naturally on $W_{\ell_j}^*$, 
we can  define a real-valued function $f_{ j}=f_{ j}(s)$ on $\Bbb R$ by
$$
f_{j}(s) \; :=\; \delta_j \,\ell_j^{-n}
\log \| \psi_{j}(t)\cdot \hat{X}_{\ell_j} \|^{}_{\operatorname{CH}(\rho^{}_{\ell_j})}, 
\qquad s \in \Bbb R.
$$
Put $\dot{f}_{j}:= df_{j}/ds$.
Since $h$ is fixed, the derivative $\dot{f}_{j}(0)$ is bounded from above by a positive constant 
 independent
of the choice of $j$ (see \cite{M}).
Hence we can define $F_{1} (\{\mu_{j}\} )\in \Bbb R \cup \{-\infty \}$ by 
$$
F_{1} (\{\mu_{j}\} )\, :=\; \lim_{s\to -\infty} \{
\varliminf_{j\to \infty}  \dot{f}_{j}(s)\},
\leqno{(2.1)}
$$
because the limit function $\varliminf_{j\to \infty}  \dot{f}_{j}(s)$  is non-decreasing  in $s$ by convexity 
of the function $f_{j}$ (cf. \cite{Zh}; see also \cite{M1}, Theorem 4.5).

\section{The sequence $\{\mu_{\gamma}\}$ generated by one test configuration}

In this section, we fix a test configuration  $\mu = (\mathcal{X},\mathcal{L}, \psi  )$ 
for $(X,L)$, of exponent $\ell$, as in the introduction.
Put $k = \ell \gamma$
for positive integers $\gamma$.
Then by the affirmative solution of equivariant Serre's conjecture, the direct image sheaves
$E^{\gamma}  :=  \pi_*\mathcal{L}^{\otimes \gamma}$,
$\gamma = 1,2,\dots$,
can be viewed as the trivial vector bundles over $\Bbb A^1$ by a $\Bbb G_m$-equivariant isomorphism 
of vector bundles 
$$
E^{\gamma} \cong \Bbb A^1\times E^{\gamma}_0,
$$
where the Hermitian metric $\rho_{k}$ on 
the fiber $E^{\gamma}_1 \; (= V_{k})$ over $1$
is chosen to be the metric as in (1.1) in the introduction, and 
the isomorphism above takes
$\rho_k$ to a Hermitian metric on the central 
fiber $E^{\gamma}_0$ which is preserved by the action of $S^1\subset \Bbb G_m$ on $E^{j}_0$ (see \cite{D2}).
Here the $\Bbb G_m$-action on $\Bbb A^1\times E^{\gamma}_0$ is induced by the $\Bbb G_m$-action on $\Bbb A^1$
and a natural $\Bbb G_m$-representation
$$
\psi_{\gamma}\,: \; \Bbb G_m \,\to  \; \operatorname{GL}(E^{\gamma}_0)\; (=\operatorname{GL}(V_{k})),
\leqno{(3.1)}
$$
where we identify $E^{\gamma}_0$ with $E^{\gamma}_1 \; (= V_{k})$ by the isomorphism above.
For $\gamma = 1$, we have $\psi_{1} = \psi$. 
Consider the sequence of
test configurations
$\mu_{\gamma }= (\mathcal{X}_{\gamma }, \mathcal{L}_{\gamma }, \psi_{\gamma})$, $\gamma = 1,2,\dots,$ 
for $(X,L)$, of exponent $k$, defined by
$$
\mathcal{X}_{\gamma } := \mathcal{X}\;\;\text{ and }\;\;\mathcal{L}_{\gamma } := \mathcal{L}^{\otimes \gamma}.
$$
Each $t\in\Bbb G_m$ viewed not as a complex number but as an element of the abstract 
group $\Bbb G_m$ inducing a biholomorphic trnsformation by $\psi_{\gamma}(t)$
on $(\mathcal{X}, \mathcal{L})$
will be denoted by $g_{\gamma}(t)$.
Let $g(t)$ denote the automorphism of $\mathcal{X}$ induced by the action of $t \in \Bbb G_m$.
Then there exists an orthonormal basis 
$\{\tau_{\alpha}\,;\, \alpha = 1,2,\dots, N_k\}$ for $V_k$ such that
$$
g_{\gamma} (t) \cdot \tau_{\alpha} \;=\; t^{-c_{\alpha}}\tau_{ \alpha},
\qquad \alpha = 1,2,\dots, N_k.
\leqno{(3.2)}
$$
where $c_{\alpha}$, $\alpha = 1,2,\dots, N_k$, are the weights of the $\Bbb G_m$-action 
on $V_k^*$. Then for the dual basis $\{\tau^*_{\alpha}\,;\, \alpha = 1,2,\dots, N_k\}$ of $V^*_k$,
we obtain
$$
g_{\gamma} (t) \cdot \tau^*_{\alpha} \;=\; t^{c_{\alpha}}\tau^*_{\alpha},
\qquad \alpha = 1,2,\dots, N_k.
$$
Each $\zeta = (\zeta^{}_{1}, \zeta^{}_{2}, \dots, \zeta^{}_{N_k})\in \Bbb C^{N_k}\setminus \{0\}$ sitting over $[\zeta ] \in \Bbb P^{N_k -1}(\Bbb C ) = \Bbb P^*(V_k )$ is viewed as  
$\Sigma_{\alpha =1}^{N_k} \zeta_{\alpha}\tau^*_{\alpha}$, so that the action  by $t \in \Bbb G_m$ 
on $\Bbb C^{N_k}$ is written as
$$
\zeta = (\zeta^{}_{1}, \zeta^{}_{2}, \dots, \zeta^{}_{N_k})\; 
\mapsto \; 
g_{\gamma}(t)\cdot \zeta = (t^{c_{1}}\zeta^{}_{1},
t^{c_{2}}\zeta^{}_{2}, \dots, 
t^{c_{N_k}}\zeta^{}_{N_k}).
\leqno{(3.3)}
$$
For the Kodaira embedding $\Phi_k : X \to \Bbb P^*(V_k )$ 
as in the introduction, we consider its image 
$X_k \,:=\, \Phi_k (X)$.
For each $p\in {\mathcal{X}}_1 \,(= X_k)$,
writing $p$ as ${\Phi}_{k} (x)$ for some $x\in X$,  we have a line $l_p$ through the origin in $\Bbb C^{N_{k}}$
associated to  the point $p = {\Phi}_{k} (x)$ in $\Bbb P^{N_{k}-1}(\Bbb C )$. 
For $t \in \Bbb R_+$, we put $p' := {g} (t)\cdot p \in {\mathcal{X}}_t$, 
where we view ${\mathcal{X}}_t$ in $ \{t\}\times \Bbb P^{N_{k}}(\Bbb C)$  
as a subvariety of $\Bbb P^{N_{k}}(\Bbb C)$ by identifying $ \{t\}\times \Bbb P^{N_{k}}(\Bbb C)$ 
with $ \Bbb P^{N_{k}}(\Bbb C)$.
Let $l_{p'}$ denote the line through the origin in $\Bbb C^{N_{k}}$ associated to
$p'$ in $\Bbb P^{N_{k}}(\Bbb C)$.
Then 
$$
l_p \to l_{p'}, \qquad {\zeta} \mapsto {g}_{\gamma} (t) \cdot {\zeta},
$$
naturally defines a map of ${\mathcal{L}}^{-\gamma}_p$ onto ${\mathcal{L}}_{p'}^{-\gamma}$.
 Thus the $\Bbb G_m$-action 
on $\mathcal{X}$ lifts to a $\Bbb G_m$-action on $\mathcal{L}^{-\gamma}$ such that
$$
\Bbb R_+ \times \mathcal{L}^{-\gamma}\to \mathcal{L}^{-\gamma}, \qquad ( t ,\zeta ) \,\mapsto  \,
{g}_{\gamma}(t)\cdot \zeta.
\leqno{(3.4)}
$$
(1) Let $\mathcal{H}_t$, $t \in \Bbb R_+$, be the set of all Hermitian metrics for the line bundle 
$\mathcal{L}_{|\mathcal{X}_t}^{-\gamma}$. 
Note that, via the identification of $\mathcal{O}_{\Bbb P^*(V_{k})}(-1)$ with $V^*_k\setminus \{0\}$,
$$
\Omega_{\operatorname{FS}}\; :=\;(n!/\ell^n)\,\Sigma_{\alpha =1}^{N_{k}} |\zeta_{\alpha}|^2
$$
defines a Hermitian metric on $\operatorname{pr}_2^*\mathcal{O}_{\Bbb P^*(V_{k})}(-1)$
and hence
on $\mathcal{L}^{-\gamma}$ when restricted to $\mathcal{X}$.
Fixing a local base for $\mathcal{L}^{\gamma}$ at $p\in X_k$, we can view each $\tau_{\alpha}$ as a
complex number.
Then ${g}(t)$ maps $ p = [\zeta ]$ to $  p' = [{g}_{\gamma}(t)\cdot \zeta ]$, and it lifts to a map sending
$q := \Omega_{\operatorname{FS}}{\,}_{|p} = (n!/\ell^n)\Sigma_{\alpha =1}^{N_{k}} |\tau_{\alpha} (p) |^2$ to 
$q' := \{{g}_{\gamma}(t) \cdot \Omega_{\operatorname{FS}}\}{\,}_{|p'}$ defined by setting
\begin{align*}
q':\, &=\,(n!/\ell^n)\Sigma_{\alpha =1}^{N_{k}}
 | {g}_{\gamma}(t)\cdot \zeta_{\alpha} |^2{\,}_{|\zeta = \tau (p)}
\,=\, (n!/\ell^n)\Sigma_{\alpha =1}^{N_{k}} t^{2c_{\alpha}} | \zeta_{\alpha} |^2{\,}_{|\zeta = \tau (p)}\\
& =\,(n!/\ell^n)\Sigma_{\alpha =1}^{N_{k}} t^{2c_{\alpha}} | \tau_{\alpha}(p) |^2.
\end{align*}
Here $\zeta = \tau (p)$ means that $\zeta_{\alpha} = \tau_{\alpha} (p)$ for all $\alpha$.
Thus $\Omega_{\operatorname{FS}}{\,}_{|X_k}\in \mathcal{H}_1$ is mapped to 
$\{{g}_{\gamma}(t) \cdot \Omega_{\operatorname{FS}}\}{\,}_{|\mathcal{X}_t} \in \mathcal{H}_t$
(cf.~\cite{M2}).

\medskip\noindent
(2) We here explain another formulation (cf.~\cite{M0}) which looks very differently but is essentially equivalent to (1) above. For each $t \in \Bbb R_+$, 
we put  $p'':= g(t^{-1})\cdot p\in \mathcal{X}_{t^{-1}}$.
Then $g(t^{-1})$ maps 
$$
p = [\zeta ] = (\tau_1 (p):\tau_2 (p):\dots :\tau_{N_k}(p))\in X_{k}
$$ 
to $p'':= [g_{\gamma}(t^{-1})\cdot \zeta ]$, 
and we have its lifting to a map taking $q := \Omega_{\operatorname{FS}}{\,}_{|p} = (n!/\ell^n)\Sigma_{\alpha =1}^{N_{k}} |\tau_{\alpha} (p) |^2$ to 
$q''= g_{\gamma}(t^{-1})\cdot q$ defined by
$$
q''  := (n!/\ell^n)\Sigma_{\alpha =1}^{N_k} | g_{\gamma}(t^{-1})\cdot \tau_{\alpha} (p)|^2 
= (n!/\ell^n)\Sigma_{\alpha =1}^{N_{k}} t^{2c_{\alpha}} | \tau_{\alpha}(p) |^2.
\leqno{(3.5)}
$$
Here, even if $t$ is replaced by $t^{-1}$ in (3.5),  $q''$  
differs from $q'$.
For $\zeta  = \Sigma_{\alpha =1}^{N_k}\zeta^{}_{\alpha}\tau^*_{\alpha}(p)$
above, $g_{\gamma}(t^{-1})\cdot \zeta = (t^{-c_{1}}\zeta^{}_{1},
t^{-c_{2}}\zeta^{}_{2}, \dots, 
t^{-c_{N_k}}\zeta^{}_{N_k})$ by (3.3), so that $q''$ in (3.5) satisfies the following:
\begin{align*}
&\langle q'',\, |g_{\gamma}(t^{-1})\cdot \zeta|^2 \rangle 
\;=\; \langle \,(n!/\ell^n) \Sigma_{\alpha =1}^{N_{k}} t^{2c_{\alpha}} | \tau_{\alpha}(p) |^2, \,
|\Sigma_{\alpha =1}^{N_k} t^{-c_{\alpha}}\zeta_{\alpha}\tau_{\alpha}^*(p)|^2\,
\rangle \\
&= \; \langle\, (n!/\ell^n) \Sigma_{\alpha =1}^{N_{k}} t^{2c_{\alpha}} | \tau_{\alpha}(p) |^2, \,
\Sigma_{\alpha=1}^{N_k}\Sigma^{N_k}_{\beta=1} t^{-c_{\alpha}-c_{\beta}} \zeta_{\alpha}\bar{\zeta}_{\beta}\tau_{\alpha}^*(p)\bar{\tau}_{\beta}^*(p)\,\rangle \\
&=\; \langle \, (n!/\ell^n) \Sigma_{\alpha =1}^{N_{k}} | \tau_{\alpha}(p) |^2, \,
\Sigma_{\alpha =1}^{N_k} |\zeta_{\alpha}|^2|\tau_{\alpha}^*(p)|^2\,\rangle \\
&=\; \langle\, (n!/\ell^n) \Sigma_{\alpha =1}^{N_{k}} | \tau_{\alpha}(p) |^2, \,
|\Sigma_{\alpha =1}^{N_k} \zeta_{\alpha}\tau_{\alpha}^*(p)|^2\,\rangle\;
=\; \langle q,\, \zeta \rangle.
\end{align*}
Hence (3.5) defines a mapping of $\mathcal{H}_1$ to $\mathcal{H}_{t^{-1}}$ naturally induced
by the map $g_{\gamma}(t^{-1}) : \mathcal{L}^{-\gamma}{}_{|X_k} \to \mathcal{L}^{-\gamma}{}_{|\mathcal{X}^{}_{t_{}^{-1}}}$
in (3.4).

\section{Behavior of $\dot{f}_{\gamma}$ for the sequence $\{\mu_{\gamma}\}$}

For a test configuration $\mu = (\mathcal{X},\mathcal{L},\psi )$ as in the previous section,
by assuming $\mu$ to be nontrivial, we here study the sequence of test configurations 
$\{\mu_{\gamma}\}$.
Put $\ln := (1/2\pi )\log$ for simplicity.
Define a Hermitian fiber norm $\phi_k $ for the real line bundle $|\mathcal{L}|^{-2/\ell}$ 
on  $\mathcal{X}$
and a K\"ahler form $\omega_k$  on $\mathcal{X}$ by
$$
\begin{cases}
&\phi_k \;:=\; (\Omega_{\operatorname{FS}})^{1/k} 
=\; \{ (n!/\ell^n) \Sigma_{\alpha =1}^{N_k}
|\zeta_{\alpha} |^2\}^{1/k} 
\\
&\omega_k \;:=\;  \sqrt{-1}\partial\bar{\partial} \ln 
\phi_k.
\end{cases}
$$
Let $\phi_k{}_{|X_k }$ denote
$\{ (n!/\ell^n) \Sigma_{\alpha =1}^{N_k}|\tau_{\alpha} |^2\}^{1/k}$, which
is justified by  $\tau_{\alpha} = \Phi_k^*\zeta_{\alpha}$.
For each $t \in \Bbb R_+$, we consider the set $\mathcal{K}_t$ of all fiber norms for the real line bundle  $|\mathcal{L}|^{-2/\ell }
\, (= |\mathcal{L}_{\gamma}|^{-2/k})$ over $\mathcal{X}_t$.
Then by (2) in Section 3,  $g_{\gamma}(t^{-1})$ maps $\mathcal{K}_1$ to 
$\mathcal{K}^{}_{t^{-1}}$. 
In particular $g_{\gamma} (t^{-1})$ maps $\phi_k{\,}_{|X_k} \in \mathcal{K}_1$  to
$$
g_{\gamma}(t^{-1})\cdot (\phi_k{}_{|X_k } )
 \; =\; \{(n!/\ell^n)\Sigma_{\alpha =1}^{N_k}\, t^{2c_{\alpha}} |\tau_{\alpha}|^2\}^{1/k},
 \leqno{(4.1)}
$$
which is viewed as an element of $\mathcal{K}^{}_{t^{-1}}$.
We now put 
$$
\bar{c}\; :=\; \Sigma_{\alpha =1}^{N_k}c_{\alpha}/N_k.
$$
Each $t \in \Bbb R_+$ viewed not as a real number but as a holomorphic transformation by 
$\psi_{\gamma}^{\operatorname{SL}}(t)$ on $V_k$
  and powers of $\mathcal{L}$
will be denoted by $\bar{g}_{\gamma}(t)$. 
For 
$$
{b}_{\alpha} := c_{\alpha} - \bar{c}, \qquad\alpha = 1,2,\dots, N_k,
\leqno{(4.2)}
$$
we have
$\bar{g}_{\gamma}(t^{-1})\cdot  \tau_{\alpha} \; =\; t^{{b}_{\alpha}}\tau_{\alpha}$,
$\alpha = 1,2,\dots, N_k$.
In view of (3.1) and (3.2) together with the definition of $\bar{c}$,  given a positive integer $r \gg 1$, we have the following asymptotic 
expansion
$$
-\,\bar{c}/ k \; = \; F_0(\mathcal{X}, \mathcal{L}) + \Sigma_{i=1}^{r} F_i (\mathcal{X}, \mathcal{L}) k^{-i}
+ O(k^{-r-1}), \qquad \gamma \gg 1.
\leqno{(4.3)}
$$
 Here for each integer $r$, 
$O(k^r)$ denotes a function $u$ satisfying $|u|\leq Ck^r$ for some positive real constant 
$C$ independent of $k$, $\alpha$ and $t$. On the other hand, by the action of $\bar{g}_{\gamma}(t)$ on 
$\mathcal{L}^{-1}$, we have
$$
\bar{g}_{\gamma}(t)\cdot (\zeta_1,\zeta_2,\dots,\zeta_{N_k})
\; =\;  (t^{b_1}\zeta_1, t^{b_2}\zeta_2, \dots , t^{b_{N_k}}\zeta_{N_k}),
$$
so that $\bar{g}_{\gamma}(t)^*\phi_k = \{(n!/\ell^n)\Sigma_{\alpha =1}^{N_k} |\bar{g}_{\gamma}(t)^*\zeta_{\alpha}|^2\}^{1/k} = \{(n!/\ell^n)\Sigma_{\alpha =1}^{N_k} t^{2b_{\alpha}}|\zeta_{\alpha}|^2\}^{1/k}$. 
Moreover by (4.2), we can write
$$
t^{-\bar{c}}\,{g}_{\gamma}(t^{-1})\cdot \tau_{\alpha}\, =\, \bar{g}_{\gamma}(t^{-1})\cdot \tau_{\alpha}\, = \,t^{b_{\alpha}}\tau_{\alpha},
\qquad \alpha = 1,2, \dots, N_k, 
$$
on $\mathcal{L}^{-1}_{\gamma}$ over $\mathcal{X}_{t^{-1}}$. Then
$t^{-2\bar{c}}{g}_{\gamma}(t^{-1})\cdot |\tau_{\alpha}|^2 
= t^{-2\bar{c}} |{g}_{\gamma}(t^{-1})\cdot \tau_{\alpha}|^2 
= t^{2b_{\alpha}}|\tau_{\alpha}|^2$, and we sum these up for $\alpha = 1,2,\dots, N_k$.
Then for the $k$-th root of the total sum,
by taking its pullback to $X_k$, we obtain
$$
\begin{cases}
&\;\;t^{-2\bar{c}/k}g(t^{-1})^*\{g_{\gamma}(t^{-1})\cdot (\phi_k{\,}_{|X_k})\}\\
&\;\;=\; \{(n!/\ell^n)\Sigma_{\alpha =1}^{N_k} t^{2b_{\alpha}}|\tau_{\alpha}|^2\}^{1/k}\; =\; \{\bar{g}_{\gamma}(t)^*\phi_k\}{\,}_{|X_k},
\end{cases}
\leqno{(4.4)}
$$
where we identify $X_k$ with $X$ via $\Phi_k$,
and $g(t^{-1})^*\{ g_{\gamma}(t^{-1})\cdot (\phi_k{\,}_{|X_k}) \}$ is viewed as a Hermitian norm on the pullback real line bundle $g(t^{-1})^*(|\mathcal{L}|^{-2/\ell}_{|\mathcal{X}_{t^{-1}}})$ over $X_k$.
Then by \cite{Ze},  
$\xi_k:=\ln \,\{ \,(\phi_k{\,}_{|X_k})/h^{-1}\}$ viewed as a function on $X$ satisfies
$$
\|\xi_k\|^{}_{C^5(X)} \;  = \;  O(k^{-2} ),
\leqno{(4.5)}
$$
where both $\phi_k{\,}_{|X_k}$ and $h^{-1}$ are viewed as Hermitian norms for the real line bundle $|L|^{-2}$ on 
$X \, (= X_k)$.
By the definition of $\xi_k$ together with (4.4),
$$
\{\bar{g}_{\gamma}(t)^*\phi_k\}{\,}_{|X_k}\, =\,\exp (2 \pi \xi_k )\cdot \kappa (t),
\leqno{(4.6)}
$$
where $\kappa (t) :=  t^{-2\bar{c}/k}\,g(t^{-1})^*\{ g_{\gamma}(t^{-1})\cdot h^{-1}\}$
and $\omega (t):= \sqrt{-1}\partial \bar{\partial}\ln \kappa (t)$.
For each $t \in \Bbb R_+$, 
we define $\theta_{t,k}$ and $\eta_{t,k}$ formally by
\begin{align*}
\;\; \theta_{t,k} &:=\; \ln\, \{\bar{g}_{\gamma}(t)^*\phi_k\}{\,}_{|X_k}\;=\;
 \ln\, (\,t^{-2\bar{c}/k}g(t^{-1})^*\{g_{\gamma}(t^{-1})\cdot (\phi_k{\,}_{|X_k})\}), 
 \tag{4.7}
 \\
\;\; \eta_{t,k}  &:= \; \ln\,\kappa (t) \qquad\;\;\;\,\quad\;\; =\;  
\ln\, (\,t^{-2\bar{c}/k}g(t^{-1})^*\{g_{\gamma}(t^{-1})\cdot h^{-1}\}),
\tag{4.8}
\end{align*}
so that 
$\{\bar{g}_{\gamma}(t)^*\omega_k\}{\,}_{|X_k} = \sqrt{-1}\partial\bar{\partial}\theta_{t,k}$ 
and $\omega (t) = \sqrt{-1}\partial\bar{\partial} \eta_{t,k}$.
In this sense, $\theta_{t,k}$ and $\eta_{t,k}$ are regarded as K\"ahler potentials 
for $\{\bar{g}_{\gamma}(t)^*\omega_k\}{\,}_{|X_k} $ and $\omega (t)$, respectively.
More precisely, in order to have functions for K\"ahler potentials, 
$t^{-2\bar{c}/k}g(t^{-1})^*\{g_{\gamma}(t^{-1})\cdot (\phi_k{\,}_{|X_k})\}$ 
and $t^{-2\bar{c}/k}g(t^{-1})^*\{g_{\gamma}(t^{-1})\cdot h^{-1}\}$
in (4.7) and (4.8) above have to be divided
by some reference Hermitian norm for $g(t^{-1})^*(|\mathcal{L}|^{-2/\ell}{}_{|\mathcal{X}^{}_{t^{-1}}})$ 
on $X_k$.
In view of the definition of $\xi_k$, subtracting (4.8) from (4.7), we obtain
$$
\theta_{t,k} - \eta_{t,k}\; =\; \ln \, \left ( \,g(t^{-1})^* \{g_{\gamma}(t^{-1})\cdot ((\phi_k{\,}_{|X_k})/ h^{-1})\}\,\right )
\; =\; \xi_k,
\leqno{(4.9)}
$$
so that $\sqrt{-1}\partial\bar{\partial}\theta_{t,k} = \omega (t) + \sqrt{-1}\partial\bar{\partial}\xi_k$.
By (4.4) and (4.7), we have $\theta_{t,k} =\, \ln\, \{ (n!/\ell^n)(\Sigma^{N_k}_{\alpha =1} \,t^{2b_{\alpha}}|\tau_{\alpha}|^2)^{1/k}\}$. We now introduce a real parameter $s$ by setting   $t = \exp (s/ \|\mu_{\gamma}\|_{\infty})$.
Then we put
$$
I_s\; := \;\delta_k\int_X ({\partial\theta_{t,k}}/{\partial s})
\, \{ (\omega (t)
+\sqrt{-1}\partial\bar{\partial}\xi_k)^n - \omega (t)^n \},
\leqno{(4.10)}
$$
where $\delta_k:=\|\mu_{\gamma}\|^{}_{\infty}/\|\mu_{\gamma}\|_1^{}$.
Put $\eta_{t}  :=  \ln\,  (\,g(t^{-1})^*\{g_{\gamma}(t^{-1})\cdot h^{-1}\})$. 
Note that the actions by $\{ g_{\gamma} (t)\,;\, t\in \Bbb G_m\}$ 
on $\mathcal{L}_{\gamma}$, $\gamma = 1,2,\dots$, are induced by the same $\Bbb G_m$-action on $\mathcal{L}$. 
Hence $\eta_t$ is independent of the choice of $\gamma$.
In view of \cite{Zh} (see also \cite{M1} and \cite{Sn}), it follows from (4.9) and (4.10)
 that
$$
\begin{cases}
&k^{-1}\dot{f}_{k} (s ) \;
=\; \delta_k \int_X ({\partial\theta_{t,k}}/{\partial s})
\, (\sqrt{-1}\partial\bar{\partial}\theta_{t,k} )^n \\
&=\,
\delta_{k}\int_X   ({\partial\theta_{t,k}}/{\partial s})
\, (\omega (t)
+\sqrt{-1}\partial\bar{\partial}\xi_k)^n\\
& = \,  I_s\,+\, \delta_k\int_X  ({\partial\theta_{t,k}}/{\partial s})
\, \omega (t)^{\,n}\, =\, I_s\,+\, \delta_k\int_X  ({\partial\eta_{t,k}}/{\partial s})
\, \omega (t)^{\,n}\\
& = \,  I_s\,  -\,   \bar{c}\,\lambda\,( k\|\mu_{\gamma}\|_{1}^{})^{-1}\,
+\,  \delta_k\int_X  ({\partial\eta_{t}}/{\partial s})
\, \omega (t)^{\,n},
\end{cases}\leqno{(4.11)}
$$
where the last equality follows from $\eta_{t,k} = \,-\,(2\bar{c}/k) \ln (t) + \eta_t$.
Since $t = \exp (s/ \|\mu_{\gamma}\|_{\infty})$, we can write
${\partial\theta_{t,k}}/{\partial s}$
in the form
$$
\frac{1}{2\pi}\,
\frac{\Sigma^{N_j}_{\alpha =1} \,2b_{\alpha} t^{2b_{\alpha}}|\tau_{\alpha}|^2}{\Sigma^{N_k}_{\alpha =1}\,  k \|\mu_{\gamma}\|_{\infty}\,
t^{2b_{\alpha}}|\tau_{\alpha}|^2}
\; =\; O(1),
$$
where $k \|\mu_{\gamma} \|_{\infty} = \max \{|b_{\alpha}|\,;\, \alpha = 1,2,\dots, N_k\}$.
Then by (4.5) and (4.10),
we obtain $k I_s = O(k^{-1})$.
Now by (4.3) and (4.11), given an integer $r \gg 1$,
$$
\begin{cases}
&\dot{f}_k (s) \,=\,  \lambda \{ \,-\,\bar{c}\, +\, k\,A(t) \,
 \}\,
\|\mu_{\gamma}\|^{\,-1}_{1} + \,kI_s \\
&=\; \lambda \{\, \tilde{F}_0(t)\,k + 
 \Sigma_{i=1}^r F_i (\mathcal{X},\mathcal{L})k^{1-i} +  O(k^{-r})
\,\}\,  \|\mu_{\gamma}\|^{\,-1}_{1} + \,kI_s,
\end{cases}
\leqno{(4.12)}
$$
where $A(t) := \lambda^{-1}\int_X t \,(\partial \eta_t/\partial t) \,\omega (t)^{\,n}$
and $\tilde{F}_0(t) := F_0(\mathcal{X},\mathcal{L}) + A(t)$. Since
$$
\omega (t) \; =\; \sqrt{-1}\partial\bar{\partial}\eta_{t,k} \; =\; \sqrt{-1}\partial\bar{\partial}\eta_{t},
$$
the function $A(t)$ is independent of the choice of $\gamma$.
Note that $k = \gamma \ell$, where $\gamma$ runs through the set of all positive integers.
Then by \cite{M},
$$
F_1(\mathcal{X},\mathcal{L}) \,+\, O(k^{-1})\; \leq \dot{f}_k (s)\;  \leq C,
\qquad  -\infty < s \leq 0,
\leqno{(4.13)}
$$
where $C$ is a positive real constant independent of the choice of $s$ and $k$.
It now follows from $\|\mu_{\gamma}\|_1 = O(1)$ that
$$
\|\mu_{\gamma}\|^{}_1 |\dot{f}_k (s)|\, =\, O(1)
\quad \text{ and } \quad \|\mu_{\gamma}\|^{}_1 k\,I_s \,=\, O(k^{-1}).
$$
Hence for each fixed $t \in \Bbb R$ with $0 < t\leq 1$, multiplying 
the equality (4.12) by $\|\mu_{\gamma}\|_1$,
and letting $k \to \infty$, 
we obtain
$$
\tilde{F}_0 (t) \; =\; 0.
\leqno{(4.14)}
$$
In particular, the function $A(t)$ takes the constant value $-F_0(\mathcal{X},\mathcal{L})$. 
Let $r$ be the smallest positive integer $i$ such that $F_i(\mathcal{X},\mathcal{L}) \neq 0$.
Here if $F_i(\mathcal{X},\mathcal{L})$ vanishes for all positive integers $i$, then 
we put $r = +\infty$ and $F_{\infty}(\mathcal{X},\mathcal{L}) =0$.
Then by (4.12) and (4.14), for all  $s\in \Bbb R$ with $s \leq 0$,
$$
\dot{f}_k (s) \; =\; \lambda \,\{ \, F_{r}(\mathcal{X},\mathcal{L}) k^{1-r}
+ O(k^{-r}) \,\}\|\mu_{\gamma}\|_1^{-1} \, +\, O(k^{-1}). 
\leqno{(4.15)}
$$
Let $r \neq +\infty$. Then by (4.13), $\|\mu_{\gamma} \|_1^{-1} \,=\,O(k^{r-1})$,
and for each fixed $s\in \Bbb R$ with $s\leq 0$, we have
$\varliminf_{\gamma \to \infty} \dot{f}_k(s)\,=\, 
\lambda\, \varliminf_{\gamma \to \infty}  
\{F_r (\mathcal{X},\mathcal{L}) k^{1-r} \|\mu_{\gamma}\|_1^{-1}\}$. Hence
$F_1 (\{\mu_{\gamma}\}) = \lim_{s\to -\infty}\varliminf_{\gamma \to \infty} \dot{f}_k(s)
= \lambda\, \varliminf_{\gamma \to \infty}  
\{F_r (\mathcal{X},\mathcal{L}) k^{1-r} \|\mu_{\gamma}\|_1^{-1}\}$. Since $k = \gamma \ell$, 
$\gamma = 1,2,\dots$, we
now obtain

\medskip\noindent
{\bf Proposition.}  
{\em \,If $\,r\neq +\infty$, then
$F_1 (\{\mu_{\gamma}\})  =  \lambda\, \varliminf_{\gamma \to \infty}  
\{F_r (\mathcal{X},\mathcal{L}) k^{1-r} \|\mu_{\gamma}\|_1^{-1}\}$
and $\|\mu_{\gamma}\|^{-1}_1 = O(\gamma^{r-1})$. If $r = +\infty$, then $F_1 (\{\mu_{\gamma}\}) = 0$.}


\medskip\noindent
{\em Remark\/  $4.16$}. It is easily seen that $F_0 (\mathcal{X},\mathcal{L})$
appearing in the asymptotic expansion (4.3) coincides with the one in the introduction.

\section{Proof of Theorem}

Let $\mu_j = (\mathcal{X}_j, \mathcal{L}_j, \psi_j )$, $ j=1,2,\dots$, be a sequence of test configurations 
for $(X,L)$ as in (1.1) in the introduction.  
In this section, we follow the arguments in the previous section.
For the time being, we fix a $j$ and let 
$$
c_{j,\alpha},\qquad \alpha = 1,2, \dots, N_{\ell_j}
$$
be the weights 
of the $\Bbb G_m$-action on $V_{\ell_j}$. 
Then there exists an orhonormal basis $\{\tau_{j,1}, \tau_{j,2},
\dots, \tau_{j,N_{\ell_j}}\}$ for $V_{\ell_j}$ such that
$g_j (t^{-1}) \cdot \tau_{\alpha} = t^{c_{j,\alpha}} \tau_{\alpha}$, $\alpha = 1,2,\dots, N_{\ell_j}$.
We now put $\bar{c}_j := c_j/N_{\ell_j}$, where
$$
{c}_j \; :=\; \Sigma_{\alpha =1}^{N_{\ell_j}} c_{j,\alpha}.
$$
Each $t \in \Bbb G_m$ viewed not as a complex number but as an element of the abstract group $\Bbb G_m$ of holomorphic transformations by $\psi_j (t)$
on $(\mathcal{X}_j, \mathcal{L}_j, V_{\ell_j})$ will be denoted by $g_j (t)$. 
Moreover, each $t \in \Bbb R_+$ viewed not as a real number but as a holomorphic transformation by $\psi^{\operatorname{SL}}_j (t)$
on $(\mathcal{X}_j,\mathcal{L}_j, V_{\ell_j})$ will be denoted by $\bar{g}_j (t)$.  By setting 
$$
b_{j,\alpha} := c_{j,\alpha} - \bar{c}_j, \qquad 
\alpha =1,2,\dots, N_{\ell_j},
$$
we have $\bar{g}_j (t^{-1} )\cdot \tau_{\alpha} = t^{b_{j,\alpha}}\tau_{\alpha}$, 
$\alpha = 1,2,\dots, N_{\ell_j}$.
By identifying $\Bbb P^*(V_{\ell_j})$ with $\Bbb P^{N_{\ell_j}-1}(\Bbb C ) = \{(z_{j,1}:z_{j,2}:\dots :z_{j,N_{\ell_j}})\}$,
we consider the Kodaira embedding 
$$
\Phi_{\ell_j}: X \hookrightarrow \Bbb P^*(V_{\ell_j}),
\quad \; x \mapsto (\tau_{j,1}(x):\tau_{j,2}(x):\dots :\tau_{j,N_{\ell_j}}(x))
$$
associated to the complete linear system $|L^{\otimes \ell_j}|$ on $X$. 
Put 
$$
\phi_j\, := \; \{\,(n!/\ell_j^n)\Sigma_{\alpha =1}^{N_{\ell_j}}\,|z_{j,\alpha}|^2\,\}^{1/\ell_j}.
$$
By identifying $X_{\ell_j}$ with $X$ via $\Phi_{\ell_j}$, 
we can view $\xi_j := \ln \{(\phi_j{}_{|X_{\ell_j}})/h^{-1})\}$ as a function on $X$ with $\|\xi_j\|_{C^5(X)} = O(\ell_j^{-2})$ (cf.~\cite{Ze}).
Define $\theta_{t,j}$ and $\eta_{t,j}$ by
$$
\begin{cases}
&\theta_{t,j}:= \; \ln ( t^{-2\bar{c}_j/\ell_j}g_j(t^{-1})^*\{g_j (t^{-1})\cdot (\phi_j {}_{|X_{\ell_j}})\}),\\
&\eta_{t,j}:= \; \ln (t^{-2\bar{c}_j/\ell_j}g_j(t^{-1})^*\{g_j (t^{-1})\cdot h^{-1}\}).
\end{cases}
$$
Put $\eta_t :=  \ln (g_j(t^{-1})^*\{g_j (t^{-1})\cdot h^{-1}\})$ and $\omega (t) := \sqrt{-1}\partial\bar{\partial}\eta_t$. For $t \in \Bbb R_+$, by introducing the parameter $s$ such that $t =  \exp (s/\|\mu_j \|_{\infty})$,
we define
$$
I_{j,s} := \delta_j \int_X (\partial \theta_{t,j}/\partial s ) 
\{(\omega (t) + \sqrt{-1}\partial\bar{\partial}\xi_j )^n - \omega (t)^n \}, \quad -\infty<s \leq 0.
$$ 
Then we can write $\theta_{t,j} = \ln \{(n!/\ell_j^n )
(\Sigma_{\alpha =1}^{N_{j_{\ell}}}t^{2b_{j,\alpha}}|\tau_{j,\alpha}|^2)^{1/\ell_j}\}$.
Put $A_j (t) := \lambda^{-1}\int_X t (\partial \eta_t/\partial t ) \omega (t)^n$. 

\medskip
Let $\pi_j : \mathcal{X}_j \to \Bbb A^1$ be the natural projection for the test configuration
$\mu_j = (\mathcal{X}_j, \mathcal{L}_j,\psi_j )$.
We now apply the argument in Section 3 to the test configuration $(\mathcal{X},\mathcal{L}, \psi ) = (\mathcal{X}_j, \mathcal{L}_j, \psi_j )$. Then by taking the direct image sheaves 
$\pi_j{}_*\mathcal{L}_j^{\otimes \gamma}$, $\gamma = 2,3,\dots$, 
we obtain a associated sequence of 
test configurations
$$
(\mathcal{X}^{}_{j;\gamma}, \mathcal{L}^{}_{j;\gamma}, \psi^{}_{j;\gamma}),
\qquad \gamma = 2,3,\dots,
$$
and for $\gamma = 1$, we put $(\mathcal{X}^{}_{j;1}, \mathcal{L}^{}_{j;1}, \psi^{}_{j;1})
= (\mathcal{X}_j, \mathcal{L}_j,\psi_j )$.
Then the arguments in Section 4 go through also for the sequence
$$
(\mathcal{X}^{}_{j;\gamma}, \mathcal{L}^{}_{j;\gamma}, \psi^{}_{j;\gamma}),
\qquad \gamma = 1,2,\dots,
$$
of test configurations.
For $\gamma =1$, by the same argument as in obtaining (4.11) and (4.12), 
$\dot{f}_j (s) \,=\, \ell_j\delta_j \int_X (\partial \theta_{t,j}/\partial s) (\sqrt{-1}\partial \bar{\partial}\theta_{t,j})^n$ is written as
$$
 \lambda \{-\bar{c}_j+ \ell_j A_j(t)\}\|\mu_j\|_1^{-1} + \ell_j I_{j,s},
$$
where by (4.14), we 
see that $A_j (t)$ is a constant function in $t \in \Bbb R_+$ with value $-F_0(\mathcal{X}_j,\mathcal{L}_j)$.
Note also that $w_{\ell_j} (\mathcal{X}_j, \mathcal{L}_j) = -{c}_j$.
Moreover by \cite{Ze}, $\|\xi_j\|_{C^5(X)} = O(\ell_j^{-2})$. Hence, letting $j\to \infty$, and then
letting $s\to -\infty$,
we now conclude that
$$
F_1 (\{\mu_j\}) \; =\; \varliminf_{j\to \infty} \frac{\lambda}{\|\mu_j\|_1}
\left \{ \frac{w_{\ell_j}(\mathcal{X}_j,\mathcal{L}_j)}{N_{\ell_j}} \,-\,\ell_j\, F_0(\mathcal{X}_j,\mathcal{L}_j)
\right \},
$$
as required.

\medskip\noindent
{\em Remark\/ $5.1$}. We here observe that the above Theorem shows that the double limit in (2.1) 
actually commutes.

\bigskip\noindent
{\footnotesize
{\sc Department of Mathematics}\newline
{\sc Osaka University} \newline
{\sc Toyonaka, Osaka, 560-0043}\newline
{\sc Japan}}

\end{document}